\documentclass[aos,preprint]{imsart}

\RequirePackage[OT1]{fontenc}
\RequirePackage{amsthm,amsmath,natbib,graphicx,array,amssymb,epic}
\RequirePackage[colorlinks]{hyperref}
\RequirePackage{hypernat}

\startlocaldefs
\numberwithin{equation}{section}
\theoremstyle{plain}
\newtheorem{thm}{Theorem}[section]
\newtheorem{definition}[thm]{Definition}
\newtheorem{proposition}[thm]{Proposition}
\newtheorem{lemma}[thm]{Lemma}

\theoremstyle{remark}
\newtheorem{remark}[thm]{Remark}
\newtheorem{example}[thm]{Example}
\endlocaldefs

\begin{document}

\begin{frontmatter}
\title{A class of statistical models to weaken independence in two-way contingency tables\protect}%\thanksref{T1}

\runtitle{Models to weaken independence}
%\thankstext{T1}{Footnote to the title with the `thankstext' command.}

\begin{aug}
\author{\fnms{Enrico}
\snm{Carlini}%\thanksref{t2}
\ead[label=e1]{enrico.carlini@polito.it}} \and
\author{\fnms{Fabio} \snm{Rapallo}\corref{}
%\thanksref{t3}
\ead[label=e2]{fabio.rapallo@mfn.unipmn.it}}
%\thankstext{t2}{First supporter of the project}
%\thankstext{t3}{Second supporter of the project}

\runauthor{E. Carlini and F. Rapallo}

\affiliation{Politecnico di Torino \and University of Eastern
Piedmont}

\address{Department of Mathematics\\
Politecnico di Torino \\
Corso Duca degli Abruzzi, 24 \\
10124 TORINO (Italy) \\
\printead{e1}\\
%\phantom{E-mail:\ }\printead*{e1}
}

\address{Department of Science and Advanced Technologies \\
University of Eastern Piedmont \\
Via Bellini, 25/g \\
15100 ALESSANDRIA (Italy) \\
\printead{e2}
}
\end{aug}

\begin{abstract}
In this paper we study a new class of statistical models for
contingency tables. We define this class of models through a
subset of the binomial equations of the classical independence
model. We use some notions from Algebraic Statistics to compute
their sufficient statistic, and to prove that they are log-linear.
Moreover, we show how to compute maximum likelihood estimates and
to perform exact inference through the Diaconis-Sturmfels
algorithm. Examples show that these models can be useful in a wide
range of applications.
\end{abstract}



\end{frontmatter}

\section{Introduction} \label{introsect}

One of the most popular statistical models for two-way contingency
tables is the independence model. It has became a reference tool
in applied research where categorical variables are concerned. In
many applications the independence model is sufficient to describe
and model the data, but this is not always the case. There are
situations where the independence model does not fit the data and
one has to detect more complex relations between the random
variables. Thus, different models have been introduced in order to
identify some structures in the contingency tables. Most of these
models belong to the class of log-linear models. Among these, we
recall the quasi-independence model, the quasi-symmetry model, the
logistic regression model. As a general reference for these models
see again \cite{agresti:02}. Such models have a wide spectrum of
applications in, e.g., biology, psychology and medicine. The books
by Fienberg \cite{fienberg:80}, Fingleton \cite{fingleton:84}, Le
\cite{le:98} and Agresti \cite{agresti:02} present a great deal of
examples with real data sets coming from the most disparate
disciplines.

A recent development in the area of statistical models for
contingency tables involves the use of some tools from Algebraic
Geometry to describe the structure and the properties of the
models. This field is currently known under the name of Algebraic
Statistics. While the first work on this direction relates to a
method for exact inference, see \cite{diaconis|sturmfels:98},
following papers have focused their attention on the geometry of
the statistical models through polynomial algebra. The algebraic
and geometric point of view in the analysis of probability models
allows us to generalize statistical models in presence of cells
with zero probability (toric models), to study its exponential
structure, and to make inference feasible also in models with
complex structure. This approach has been particularly useful in
the fields of log-linear and graphical models. Some relevant works
on these recent topics are \cite{geiger|meek|sturmfels:06},
\cite{garcia|stillman|sturmfels:05},
\cite{geiger|heckerman|king|meek:01}, and \cite{rapallo:07}. An
exposition of such theory, with a view toward applications to
computational biology, can be found in
\cite{patcher|sturmfels:05}.

The theoretical advances mentioned above also have a computational
counterpart. In fact, many symbolic softwares traditionally
conceived for polynomial algebra now include special functions or
packages specifically designed for Algebraic Statistics, see e.g.
{C}o{C}o{A} \cite{cocoa}, 4ti2 \cite{4ti2}, and {L}att{E}
\cite{latte}.

In this paper we consider statistical models for two-way
contingency tables with strictly positive cell probabilities. We
introduce a class of models in order to weaken independence,
starting from the binomial representation of the independence
model. The independence statement means that the table of
probabilities has rank $1$, and therefore that all $2 \times 2$
minors vanish. In the strictly positive case, this is equivalent
to the vanishing of all $2 \times 2$ adjacent minors. Our models,
which we call {\em weakened independence models}, are defined
through a subset of the independence binomial equations. As a
consequence, the independence statements hold locally and the
resulting models allow us to identify local patterns of
independence in contingency tables. We study the main properties
of such models. In particular, we prove that they belong to the
class of log-linear models, and we determine their sufficient
statistic. Moreover, we compute the corresponding Markov bases, in
order to apply the Diaconis-Sturmfels algorithm without symbolic
computations. The relevance of our theory is emphasized by some
examples on real data sets. We also show that our models have
connections with a problem recently stated by Bernd Sturmfels in
the field of probability models for Computational Biology, the
so-called ``$100$ Swiss Francs Problem'', see \cite{sturmfels:07}.

While most of the papers in Algebraic Statistics uses algebraic
and geometric methods to describe and analyze existing statistical
models, or to make exact inference, the main focus of this paper
is the definition of a new class of models, by exploiting the
Algebraic Statistics way of thinking.

Notice that we restrict the analysis to adjacent minors.
Therefore, the applications are mainly concerned with binary or
ordinal random variables. At the end of the paper we will give
some pointers to follow-ups and extensions of this work.

In Section \ref{defsect} we define the weakened independence
models and we give some examples, while in Section \ref{suffsect}
we provide the computation of a sufficient statistic. In Section
\ref{expsect} we prove that these models belong to the class of
toric models (and therefore they are log-linear for strictly
positive probabilities), and we explicitly write down some
consequences, such as a canonical parametrization of the models.
In Section \ref{markovsect} we compute the Markov bases for
weakened independence models and we present some examples with
real data. In particular, Example \ref{swiss} is devoted to the
discussion of some interesting relationships between our models
and the ``$100$ Swiss Francs Problem''. Section \ref{finremsect}
highlights the main contributions of our theory and provides some
pointers to future developments.

\section{Definitions} \label{defsect}

A two-way contingency table collects data from a sample where two
categorical variables, say $X$ and $Y$, are measured. Suppose that
$X$ has $I$ levels and $Y$ has $J$ levels. The sample space for a
sample of size one is ${\mathcal X}=\{1, \ldots ,I\}\times \{1
,\ldots, J\}$ and a joint probability distribution for an $I
\times J$ contingency table is a table of raw probabilities
$(p_{i,j})_{i=1,\ldots,I,j=1,\ldots,J}$ in the simplex
\begin{equation*}
\Delta = \left\{(p_{i,j})_{i=1,\ldots,I,j=1,\ldots,J} \in {\mathbb
R}^{I \times J}_+ \ : \ \sum_{i,j} p_{i,j} = 1 \right\} \, .
\end{equation*}
A statistical model for an $I \times J$ contingency table is then
a subset of $\Delta$ defined through equations on the raw
probabilities $p_{1,1}, \ldots, p_{I,J}$. In this paper, we do not
allow any $p_{i,j}$ to be zero, and we assume strict positivity of
all probabilities.

The independence model can be defined in parametric form through
the power product representation, i.e. by the set of equations
\begin{equation} \label{indep1}
p_{i,j} = \zeta_0 \zeta_{i}^X \zeta_{j}^Y
\end{equation}
for $i=1, \ldots, I$ and $j=1, \ldots, J$, where $\zeta_{i}^X$ and
$\zeta_{j}^Y$ are unrestricted positive parameters and $\zeta_0$
is the normalizing constant, see
\cite{pistone|riccomagno|wynn:01}. In term of log-probabilities,
Eq. \eqref{indep1} assumes the most familiar form
\begin{equation} \label{indeplog}
\log p_{i,j} = \lambda + \lambda_i^X + \lambda_j^Y
\end{equation}
where $\lambda=\log\zeta_{0}$, $\lambda_i^X=\log\zeta_{i}^X$ for
$i=1, \ldots, I$ and $\lambda_j^Y=\log\zeta_{j}^Y$ for $j=1,
\ldots J$. As an equivalent representation, one can derive
implicit formulae on the raw probabilities $p_{i,j}$. Eliminating
the $\zeta$ variables from Eq. \eqref{indep1}, one obtains the set
of equations below:
\begin{equation} \label{indep2}
p_{i,j}p_{k,m} - p_{i,m}p_{k,j} = 0
\end{equation}
for all $1 \leq i < k \leq I$ and $1 \leq j < m \leq J$. In other
words, in the independence model all $2 \times 2$ minors of the
table vanish. It is well known, see e.g. \cite{agresti:02}, that
in the positive case, the equalities in Eq. \eqref{indep2} are
redundant and it is enough to set to zero the adjacent minors:
\begin{equation} \label{indep3}
p_{i,j}p_{i+1,j+1} - p_{i+1,j}p_{i,j+1}
\end{equation}
for all $1 \leq i < I$ and $1 \leq j < J$.

\begin{remark}
In the framework of toric models as defined in
\cite{pistone|riccomagno|wynn:01}, where structural zeros are
allowed, the implicit representations \eqref{indep2} and
\eqref{indep3} are not equivalent, as they differ on the boundary.
For a description of such phenomenon, see \cite{rapallo:07}.
\end{remark}

In algebraic terms, let
\begin{equation*}
{\mathcal C} = \{p_{i,j}p_{i+1,j+1} - p_{i+1,j}p_{i,j+1} \ : \ 1
\leq i < I, \ 1 \leq j < J \} \, .
\end{equation*}
The set ${\mathcal C}$ is the set of all $2 \times 2$ adjacent
minors of the table of probabilities. Moreover, let ${\mathbb
R}[p]$ be the polynomial ring in $I \times J$ indeterminates with
real coefficients.

From the geometric point of view, the independence model is the
variety
\begin{equation*}
V_{\mathcal C} = \{p_{i,j} \ : \ {\mathcal C} = 0 \} \cap \Delta
\, ,
\end{equation*}
i.e., the set of the points of the simplex where all binomials in
${\mathcal C}$ vanish.

The choice of a subset of ${\mathcal C}$ leads us to the
definition of a new class of models.

\begin{definition} \label{weakdef}
Let ${\mathcal B}$ be a subset of ${\mathcal C}$. The {\em
${\mathcal B}$-weakened independence model} is the variety
\begin{equation*}
V_{\mathcal B} = \{p_{i,j} \ : \ {\mathcal B} = 0 \} \cap \Delta
\, .
\end{equation*}
\end{definition}

Of course, $V_{\mathcal C} \subseteq V_{\mathcal B}$ for all
subsets ${\mathcal B}$ of ${\mathcal C}$. The meaning of the class
of models in Definition \ref{weakdef} is quite simple. In fact,
the choice of a given set of minors means that we allow the
binomial independence statements to hold locally, i.e., we
determine patterns of independence.

\begin{example}
As a first applications, we consider a $2 \times J$ contingency
table. A table of this kind could derive, e.g., from the
observation of a binary random variable $X$ at different times.

The model defined through the set of binomials
\begin{equation*}
{\mathcal B} = \{p_{1,1}p_{2,2}-p_{1,2}p_{2,1},
p_{1,2}p_{2,3}-p_{1,3}p_{2,2}, \ldots , p_{1, j'-1}p_{2,j'}-p_{1,
j'}p_{2,j'-1} \} \, ,
\end{equation*}
where $j'<J$, is presented in Figure \ref{figlogreg}. This choice
of ${\mathcal B}$ means that there is independence between $X$ and
the time up to the instant $j'$ and not after.
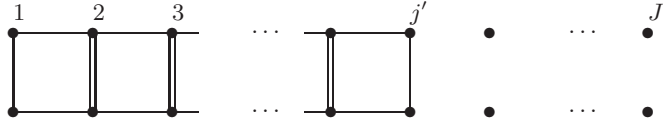
\begin{figure}
\begin{center}
\begin{picture}(260,70)(0,0)
\put(10,10){\circle*{4}}\put(10,10){\line(1,0){30}}
\put(10,10){\line(0,1){30}} \put(10,45){$1$}

\put(10,40){\circle*{4}}\put(10,40){\line(1,0){30}}

\put(40,10){\circle*{4}}\put(40,10){\line(1,0){30}}
\put(39,10){\line(0,1){30}}\put(41,10){\line(0,1){30}}
\put(40,45){$2$}

\put(40,40){\circle*{4}}\put(40,40){\line(1,0){30}}

\put(70,10){\circle*{4}}\put(70,10){\line(1,0){10}}
\put(69,10){\line(0,1){30}}\put(71,10){\line(0,1){30}}
\put(70,45){$3$}

\put(70,40){\circle*{4}}\put(70,40){\line(1,0){10}}

\put(100,10){$\ldots$} \put(120,10){\line(1,0){10}}

\put(100,40){$\ldots$} \put(120,40){\line(1,0){10}}

\put(130,10){\circle*{4}}\put(130,10){\line(1,0){30}}
\put(129,10){\line(0,1){30}}\put(131,10){\line(0,1){30}}

\put(130,40){\circle*{4}}\put(130,40){\line(1,0){30}}

\put(160,10){\circle*{4}}
\put(160,10){\line(0,1){30}}\put(160,45){$j'$}

\put(160,40){\circle*{4}}

\put(190,10){\circle*{4}}

\put(190,40){\circle*{4}}

\put(220,10){$\ldots$}

\put(220,40){$\ldots$}

\put(250,10){\circle*{4}}

\put(250,40){\circle*{4}} \put(250,45){$J$}

\end{picture}
\end{center}
\caption[]{Binomials for a change-point problem in logistic
regression.}\label{figlogreg}
\end{figure}
In literature, the point $j'$ in this model refers to the
detection of the change-point in a logistic regression model. A
recent paper about this topic is \cite{gurevich|vexler:05}.
\end{example}

\begin{example} \label{fivebyfive}
Let us consider a $I \times I$ contingency table. A table of this
kind could derive from a rater agreement analysis. Suppose that
$2$ raters independently classify $n$ objects using a nominal or
ordinal scale with $I$ categories. If we set
\begin{equation*}
{\mathcal B} = \{p_{1,1}p_{2,2}-p_{1,2}p_{2,1} \}
\end{equation*}
the corresponding model yields that categories $1$ and $2$ are
indistinguishable. A reference for the notion of category
indistinguishability is, e.g., \cite{darroch|mccloud:86}. This
model can be generalized using the set of binomials
\begin{equation*}
{\mathcal B} = \{p_{i,j}p_{i+1,j+1}-p_{i,j+1}p_{i+1,j} \ : \ 1
\leq i \leq i', \ 1 \leq j \leq i' \}
\end{equation*}
meaning that the categories $1, \ldots, i'$ are indistinguishable.
The first paper in the direction of modelling patterns of
agreement is \cite{agresti:92}. An example with $5$ categories and
$3$ undistinguishable categories is presented in Figure
\ref{fig5x5}.
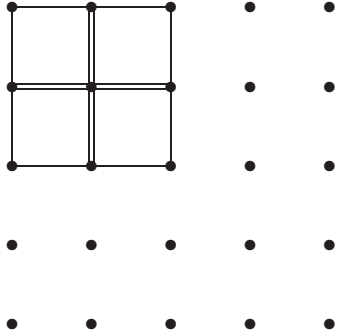
\begin{figure}
\begin{center}
\begin{picture}(140,140)(0,0)
\put(10,10){\circle*{4}}

\put(10,40){\circle*{4}}

\put(10,70){\circle*{4}} \put(10,70){\line(0,1){30}}
\put(10,70){\line(1,0){30}}

\put(10,100){\circle*{4}} \put(10,100){\line(0,1){30}}
\put(10,99){\line(1,0){30}} \put(10,101){\line(1,0){30}}

\put(10,130){\circle*{4}} \put(10,130){\line(1,0){30}}

\put(40,10){\circle*{4}}

\put(40,40){\circle*{4}}

\put(40,70){\circle*{4}} \put(39,70){\line(0,1){30}}
\put(41,70){\line(0,1){30}} \put(40,70){\line(1,0){30}}

\put(40,100){\circle*{4}} \put(39,100){\line(0,1){30}}
\put(41,100){\line(0,1){30}} \put(40,99){\line(1,0){30}}
\put(40,101){\line(1,0){30}}

\put(40,130){\circle*{4}} \put(40,130){\line(1,0){30}}

\put(70,10){\circle*{4}}

\put(70,40){\circle*{4}}

\put(70,70){\circle*{4}} \put(70,70){\line(0,1){30}}

\put(70,100){\circle*{4}} \put(70,100){\line(0,1){30}}

\put(70,130){\circle*{4}}

\put(100,10){\circle*{4}}

\put(100,40){\circle*{4}}

\put(100,70){\circle*{4}}

\put(100,100){\circle*{4}}

\put(100,130){\circle*{4}}

\put(130,10){\circle*{4}}

\put(130,40){\circle*{4}}

\put(130,70){\circle*{4}}

\put(130,100){\circle*{4}}

\put(130,130){\circle*{4}}
\end{picture}
\end{center}
\caption[]{Binomials for Example \ref{fivebyfive}.}\label{fig5x5}
\end{figure}
More examples on the models for rater agreement problems will be
presented later in the paper.
\end{example}

\begin{remark}
In the next sections, our approach will proceed somehow backwards
with respect to the classical log-linear models theory. In fact,
we will define the model through the binomials and then we will
use them to determine a sufficient statistic and a
parametrization.
\end{remark}

\section{Sufficient statistic} \label{suffsect}

As noticed in the Introduction, the independence model is defined
through the log-linear form in Eq. \eqref{indeplog}. One can
easily check that for the independence model a sufficient
statistic $T$ for the sample of size $1$ is given by the indicator
functions of the $I$ rows and the indicator functions of the $J$
columns. More precisely, we denote the indicator function of the
$i$-th row by ${\mathbb I}_{(i,+)}$ and the indicator function of
the $j$-th column by ${\mathbb I}_{(+,j)}$. Writing the sample
space as ${\mathcal X}= \{ 1, \ldots , I\} \times \{1, \ldots , J
\}$, a sufficient statistic for the independence model is
\begin{equation*}
T=\left({\mathbb I}_{(1,+)} , \ldots , {\mathbb I}_{(I,+)},
{\mathbb I}_{(+,1)}, \ldots, {\mathbb I}_{(+,J)} \right) \, .
\end{equation*}
A single observation is an element of the sample space ${\mathcal
X}$ and its table has a single count of $1$ in one cell and $0$
otherwise. This observation yields a value of $1$ in the
corresponding row and column indicator functions in $T$.

Therefore, the sufficient statistic $T$ for a sample of size $1$
is a linear map from ${\mathcal X}$ to ${\mathbb N}^{I+J}$. The
function $T$ can be extended to a linear homomorphism $T: {\mathbb
R}^{IJ} \rightarrow {\mathbb R}^{I + J}$.

In Section \ref{expsect} we will prove that weakened independence
models, as the independence model, are log-linear. Thus, the
sufficient statistic for a sample of size $n$ is the sum of the
sufficient statistics of all components of the sample and it will
be formed by the sum of appropriate cell counts, as familiar in
the field of categorical data analysis, see e.g.
\cite{agresti:02}. However, in this section it is more convenient
to work with a sample of size one and with the indicator
functions. This approach has been fruitfully used in
\cite{haberman:74} and, more recently, in
\cite{pistone|riccomagno|wynn:01}.

Hereinafter, we write the table as a column vector, i.e. the table
of probabilities is written as
\begin{equation*}
p=\left( p_{1,1}, \ldots, p_{1,J}, \ldots , p_{I,1}, \ldots,
p_{I,J} \right)^t \, ,
\end{equation*}
where $t$ denotes the transposition. Moreover, we use a vector
notation, i.e. we write a binomial in the form $p^a-p^b$, meaning
$p_{1,1}^{a_{1,1}} \cdots p_{I,J}^{a_{I,J}} - p_{1,1}^{b_{1,1}}
\cdots p_{I,J}^{b_{I,J}}$.

We briefly review the relationship between the sufficient
statistic and the binomials in Eq. \eqref{indep3}. Writing the
table as a column vector of length $IJ$, the matrix representation
of $T$ is a matrix $A_{\mathcal C}$. This matrix has size $IJ
\times (I + J)$ and its rank is $I+J-1$.

Moreover, consider the log-vector of a $2 \times 2$ minor to be
defined in the following way:
\begin{eqnarray*}
\Lambda  : &  {\mathbb R}[p] & \longrightarrow {\mathbb R}^{IJ} \\
& p^a - p^b & \longmapsto a-b
\end{eqnarray*}

We denote by $Z_{\mathcal C}$ the sub-vector space of ${\mathbb
R}^{IJ}$ generated by the vectors $\Lambda(m)$, for all $2 \times
2$ adjacent minors $m$. It is well known, see for example
\cite{bishop|fienberg|holland:75}, that $Z_{\mathcal C}$ has
dimension $(I-1)(J-1)$ and the sequence of log-vectors
$\Lambda(m)$ with $m \in {\mathcal C}$ is a sequence of
$(I-1)(J-1)$ linearly independent vectors orthogonal to
$A_{\mathcal C}$. Hence the column space $A_{\mathcal C}$ is the
orthogonal of $Z_{\mathcal C}$ Thus, from a vector-space
perspective, the exponents of the binomials are the orthogonal
complement of the matrix $A_{\mathcal C}$. In the sequel, we will
use the same symbol to denote a matrix $A$ and the sub-vector
space of ${\mathbb R}^{IJ}$ generated by the columns of $A$,
although this should be considered as a slight abuse of notation.

The procedure described above is quite general and it provides a
method to actually compute the relevant binomials of a statistical
model with a given sufficient statistic. For more details, see
\cite{pistone|riccomagno|wynn:01}.

In order to analyze the weakened independence models in Definition
\ref{weakdef}, we use the theory sketched above for the
independence model. We start with a set of binomials, we compute a
sufficient statistic and the parametric representation of the
model.

\begin{remark}
We will prove in Section \ref{expsect} that weakened independence
models are log-linear. Therefore, the orthogonal to the
log-vectors of the chosen binomials is the matrix representation
of a sufficient statistic. In order to keep notation as simple as
possible, we call this orthogonal a sufficient statistic even
before showing that the models are log-linear.
\end{remark}

\begin{lemma} \label{logveclemma}
The log-vectors of $d$ distinct adjacent minors are linearly
independent.
\end{lemma}
\begin{proof}
Let ${\mathcal B}_d$ be a set of $d$ distinct adjacent minors and
let $L_d$ be the set of their log-vectors. We proceed by induction
on $d$. For $d=1$ the statement is clearly true. We assume that
the elements of $L_i$ are linearly independent for all $i<d$ and
we will show that the same holds for $L_d$. Let $m\in {\mathcal
B}_d$ be the minor involving the indeterminate having the
lex-smallest index, say $p_{\bar i,\bar j}$, and notice that no
other element in ${\mathcal B}_d$ involves $p_{\bar i,\bar j}$.
Let $l=\Lambda(m)\in L_d$ and notice that $l$ is not a linear
combination of the element of $L_d\setminus\{l\}$ which are
linearly independent by hypothesis. Hence the element of $L_d$ are
linearly independent.
\end{proof}

\begin{remark}
Lemma \ref{logveclemma} is false when we consider log-vectors of
non-adjacent minors. As a counterexample, take a $2 \times 3$
table and all three minors.
\end{remark}

Now, consider a ${\mathcal B}$-weakened independence model with
set of adjacent minors ${\mathcal B}$ of cardinality $m$. Let
$Z_{\mathcal B}$ be the matrix of the log-vectors of the adjacent
minors in ${\mathcal B}$. In view of Lemma \ref{logveclemma}, the
orthogonal of $Z_{\mathcal B}$ has dimension $(IJ-m)$. Thus, the
explicit computation of $A_{\mathcal B}$, the orthogonal of
$Z_{\mathcal B}$, requires to find at least $(IJ-m)$ vectors
orthogonal to $Z_{\mathcal B}$. Although this can be done simply
with a linear algebra algorithm, it is very useful to investigate
the structure of the the matrix $A_{\mathcal B}$.

Given a ${\mathcal B}$-weakened independence model for $I \times
J$ contingency tables, we define a graph in the following way.

\begin{definition}
Given a set ${\mathcal B}$ of adjacent minors, we define a graph
$G_{\mathcal B}$ as follows: the set of vertices is the set of
cells and each binomial defines 4 edges. The binomial
$p_{i,j}p_{i+1,j+1}-p_{i,j+1}p_{i+1,j}$ defines the edges
$(i,j)\leftrightarrow(i+1,j)$, $(i+1,j)\leftrightarrow(i+1,j+1)$,
$(i,j+1)\leftrightarrow(i+1,j+1)$ and
$(i,j)\leftrightarrow(i,j+1)$.
\end{definition}

The edges associated to a binomial are the $4$ sides of the square
with vertices on the 4 cells involved in the binomial.

\begin{definition}
A cell $(i,j)$ is a {\em free cell} if no edge of $G_{\mathcal B}$
involves ${(i,j)}$.
\end{definition}

Equivalently, a cell $(i,j)$ is a free cell if and only if the
indeterminate $p_{i,j}$ does not appear in any of the binomials in
${\mathcal B}$.

\begin{definition}
The sequence of cells $(i,j), (i,j+1), \ldots (i,j+h)$ is a {\em
connected component} of the $i$-th row if each pair of consecutive
cells is connected by an edge of $G_{\mathcal B}$. The sequence
forms a {\em maximal connected row component} ($MCR$) if the
sequence is no more connected when one adds $(i,j-1)$ or
$(i,j+h+1)$.
\end{definition}

One can define similarly the maximal connected column component
($MCC$). We illustrate the definitions above with an example.

\begin{example} \label{ex4per4}
In the model for a $4 \times 4$ contingency table defined through
the binomials in Figure \ref{fig4x4}, we have $4$ $MCR$s, $5$
$MCC$s and 2 free cells.
\begin{figure}
\begin{center}
\begin{picture}(110,110)(0,0)
\put(10,10){\circle*{4}}

\put(10,40){\circle*{4}}

\put(10,70){\circle*{4}} \put(10,70){\line(0,1){30}}
\put(10,70){\line(1,0){30}}

\put(10,100){\circle*{4}} \put(10,100){\line(1,0){30}}

\put(40,10){\circle*{4}} \put(40,10){\line(0,1){30}}
\put(40,10){\line(1,0){30}}

\put(40,40){\circle*{4}} \put(40,40){\line(0,1){30}}
\put(40,39){\line(1,0){30}} \put(40,41){\line(1,0){30}}

\put(40,70){\circle*{4}} \put(39,70){\line(0,1){30}}
\put(41,70){\line(0,1){30}} \put(40,69){\line(1,0){30}}
\put(40,71){\line(1,0){30}}

\put(40,100){\circle*{4}} \put(40,100){\line(1,0){30}}

\put(70,10){\circle*{4}}  \put(69,10){\line(0,1){30}}
\put(71,10){\line(0,1){30}} \put(70,10){\line(1,0){30}}

\put(70,40){\circle*{4}} \put(70,40){\line(0,1){30}}
\put(70,40){\line(1,0){30}}

\put(70,70){\circle*{4}} \put(69,70){\line(0,1){30}}
\put(71,70){\line(0,1){30}} \put(70,70){\line(1,0){30}}

\put(70,100){\circle*{4}} \put(70,100){\line(1,0){30}}

\put(100,10){\circle*{4}} \put(100,10){\line(0,1){30}}

\put(100,40){\circle*{4}}

\put(100,70){\circle*{4}} \put(100,70){\line(0,1){30}}

\put(100,100){\circle*{4}}
\end{picture}
\end{center}
\caption[]{Binomials for Example \ref{ex4per4}.}\label{fig4x4}
\end{figure}
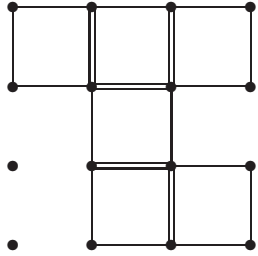
\end{example}

\begin{proposition} \label{orthoprop}
Consider a ${\mathcal B}$-weakened independence model with set of
binomials ${\mathcal B}$ and let $Z_{\mathcal B}$ be the matrix of
the log-vectors of the minors in ${\mathcal B}$. The indicator
vectors of the free cells, the indicator vectors of the $MCR$s and
the indicator vectors of the $MCC$s are orthogonal to the column
space $Z_{\mathcal B}$.
\end{proposition}
\begin{proof}
If $(i,j)$ is a free cell, then no monomial in ${\mathcal B}$
involves the corresponding variable. Hence the indicator vector of
$(i,j)$ is orthogonal to the column space of $Z_{\mathcal B}$.
Given a $MCR$, its indicator function is clearly orthogonal to the
columns of $Z_{\mathcal B}$ corresponding to minors not involving
the cells of the $MCR$. If a minor involves a cell of the $MCR$,
then it involves two cells with alternating signs. A similar
argument works for $MCC$s. Hence the orthogonality follows.
\end{proof}

Now, two questions arise: one about the linear independence of the
vectors defined in Proposition \ref{orthoprop} and the other about
the dimension of the sub-vector space generated by such vectors.
In other words, we have to investigate whether these vectors
generate the space orthogonal to $Z_{\mathcal B}$ or not. Let us
start with two simple examples.

\begin{example} \label{exdiagminors}
Consider a weakened independence model for $4 \times 4$ tables
defined through the adjacent minors in Figure \ref{fig3minors}.
\begin{figure}
\begin{center}
\begin{picture}(110,110)(0,0)
\put(10,10){\circle*{4}}

\put(10,40){\circle*{4}}

\put(10,70){\circle*{4}} \put(10,70){\line(0,1){30}}
\put(10,70){\line(1,0){30}}

\put(10,100){\circle*{4}} \put(10,100){\line(1,0){30}}

\put(40,10){\circle*{4}}

\put(40,40){\circle*{4}} \put(40,40){\line(0,1){30}}
\put(40,40){\line(1,0){30}}

\put(40,70){\circle*{4}} \put(40,70){\line(1,0){30}}
\put(40,70){\line(0,1){30}}

\put(40,100){\circle*{4}}

\put(70,10){\circle*{4}}  \put(70,10){\line(0,1){30}}
\put(70,10){\line(1,0){30}}

\put(70,40){\circle*{4}} \put(70,40){\line(0,1){30}}
\put(70,40){\line(1,0){30}}

\put(70,70){\circle*{4}}

\put(70,100){\circle*{4}}

\put(100,10){\circle*{4}} \put(100,10){\line(0,1){30}}

\put(100,40){\circle*{4}}

\put(100,70){\circle*{4}}

\put(100,100){\circle*{4}}
\end{picture}
\end{center}
\caption[]{Binomials for Example
\ref{exdiagminors}.}\label{fig3minors}
\end{figure}
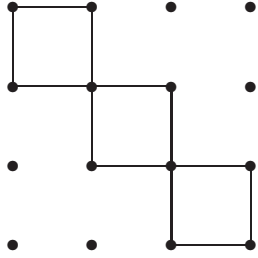
In this situation, $Z_{\mathcal B}$ has rank 3 and there are $4$
$MCR$s, $4$ $MCC$s and $6$ free cells. Here, the $14$ vectors
corresponding to the $MCR$s, to the $MCC$s and to the free cells
generate a sub-vector space of dimension $13$. Thus, they are
enough to define the matrix $A_{\mathcal B}$.
\end{example}

\begin{example} \label{patexample}
Consider now a weakened independence model for $4 \times 4$ tables
defined through the adjacent minors in Figure \ref{figbuco}.
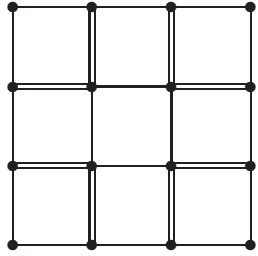
\begin{figure}
\begin{center}
\begin{picture}(110,110)(0,0)
\put(10,10){\circle*{4}} \put(10,10){\line(0,1){30}}
\put(10,10){\line(1,0){30}}

\put(10,40){\circle*{4}} \put(10,40){\line(0,1){30}}
\put(10,39){\line(1,0){30}} \put(10,41){\line(1,0){30}}

\put(10,70){\circle*{4}} \put(10,70){\line(0,1){30}}
\put(10,69){\line(1,0){30}} \put(10,71){\line(1,0){30}}

\put(10,100){\circle*{4}} \put(10,100){\line(1,0){30}}

\put(40,10){\circle*{4}} \put(39,10){\line(0,1){30}}
\put(41,10){\line(0,1){30}} \put(40,10){\line(1,0){30}}

\put(40,40){\circle*{4}} \put(40,40){\line(0,1){30}}
\put(40,40){\line(1,0){30}}

\put(40,70){\circle*{4}} \put(39,70){\line(0,1){30}}
\put(41,70){\line(0,1){30}} \put(40,70){\line(1,0){30}}

\put(40,100){\circle*{4}} \put(40,100){\line(1,0){30}}

\put(70,10){\circle*{4}}  \put(69,10){\line(0,1){30}}
\put(71,10){\line(0,1){30}} \put(70,10){\line(1,0){30}}

\put(70,40){\circle*{4}} \put(70,40){\line(0,1){30}}
\put(70,39){\line(1,0){30}} \put(70,41){\line(1,0){30}}

\put(70,70){\circle*{4}}  \put(69,70){\line(0,1){30}}
\put(71,70){\line(0,1){30}} \put(70,69){\line(1,0){30}}
\put(70,71){\line(1,0){30}}

\put(70,100){\circle*{4}} \put(70,100){\line(1,0){30}}

\put(100,10){\circle*{4}} \put(100,10){\line(0,1){30}}

\put(100,40){\circle*{4}} \put(100,40){\line(0,1){30}}

\put(100,70){\circle*{4}} \put(100,70){\line(0,1){30}}

\put(100,100){\circle*{4}}
\end{picture}
\end{center}
\caption[]{Binomials for Example \ref{patexample}.}
\label{figbuco}
\end{figure}
The model above only leaves out one minor, namely the central one.
In this case, $Z_{\mathcal B}$ has rank 8 and there are $4$ $MCR$s
and 4 $MCC$s. The $8$ vectors corresponding to the $MCR$s and to
the $MCC$s generate a sub-vector space of dimension $7$ and
therefore they are not enough to generate the orthogonal space
$A_{\mathcal B}$.
\end{example}

The dimension of the vector space generated by the indicator
function of the $MCR$s, the $MCC$s and the free cells can be
computed and we have the following results.

\begin{proposition} \label{linindprop}
For any connected component of ${\mathcal B}$ with $r$ $MCR$s and
$c$ $MCC$s, the vector space generated by the $MCR$s and by the
$MCC$s has dimension $(r+c-1)$.
\end{proposition}
\begin{proof}
Clearly the log-vectors of the $MCR$s and of the $MCC$s are not
linearly independent as their sums are equal. To show that this is
the only relation we proceed by induction on the number of minors
in the connected component. Let this number be $d$. If $d=1$ the
result is trivial. Now, assume that the result holds for $d$. If
the connected component involves $d+1$ minors, let
$r_1,\ldots,r_t$ be the indicator functions of the $MCR$s and
$c_1,\ldots,c_s$ be the indicator functions of the $MCC$s. Also
assume that $c_1$ and $r_1$ involve the lex-smallest cell. Notice
that $c_1$ and $r_1$ are the only vectors involving this cell.
Given a linear combination
\begin{equation} \label{lincomb}
\sum \lambda_i c_i=\sum \mu_i r_i
\end{equation}
we must have $\lambda_1=\mu_1$. Then the linear combination
\eqref{lincomb}  can be read in ${\mathcal B}'={\mathcal
B}\setminus\{\overline m\}$, where $\overline m$ is the minor
involving the lex-smallest cell and the $c_i$'s and the $r_i$'s
represent the log-vectors of the $MCR$s and the $MCC$s of
${\mathcal B'}$. By the inductive hypothesis we get
$\lambda_i=\mu_i=1$ for all $i$.
\end{proof}

As distinct connected components and free cells act on spaces
which are orthogonal to each other, Proposition \ref{linindprop}
leads to the following corollary.

\begin{thm} \label{maincounter}
Consider a ${\mathcal B}$-weakened independence model defined by a
set of binomials ${\mathcal B}$ whose graph has $k$ connected
components and with $r$ $MCR$s, $c$ $MCC$s and $f$ free cells. The
dimension of the vector space generated by $MCR$s, $MCC$s and the
indicator functions of the free cells is $(r+c+f-k)$.
\end{thm}
\begin{proof}
The indicators of the free cells are clearly independent with the
indicators of the $MCC$s and of the $MRC$s. Moreover, indicators
of $MCC$s and $MCR$s of different connected components are
linearly independent as they do not share any cell. By Proposition
\ref{linindprop} each connected component gives exactly one
relation among the indicators of the $MCR$s and the $MCC$s. Hence
the result follows.
\end{proof}

In the results above, we have addressed dimensional issues. Now,
we use them to find a procedure to determine a sufficient
statistic. Moreover, the examples of this section show that in
some cases the vectors of $MCC$s, $MCR$s and free cells are
sufficient to generate the space orthogonal to $Z_{\mathcal B}$.
Clearly, these vectors are not sufficient when the graph
$G_{\mathcal B}$ of the binomials in ${\mathcal B}$ present a
hole, i.e., when we remove some minors with 4 double edges from
the complete set of binomials ${\mathcal C}$. Removing such a
minor adds a new vector to the orthogonal. On the other hand, it
does not add anything in terms of $MCR$s, $MCC$s and free cells.

Thus, the last part of this section is devoted to actually find a
sufficient statistic for a generic weakened independence model.
The key idea is to start from the complete set of adjacent minors
${\mathcal C}$ and to remove minors iteratively. This approach is
motivated by the fact that for the complete set ${\mathcal C}$ a
sufficient statistic is known to be formed by the row sums and the
column sums, as extensively discussed in Section \ref{defsect}.

We begin our analysis from a simple case. Namely, we consider a
set of binomials ${\mathcal B}$ with given sufficient statistic
$A_{{\mathcal B}}$ and we investigate the behavior of the
sufficient statistic when we remove one minor $m$ from ${\mathcal
B}$, i.e., when the set of binomials is ${\mathcal B}'={\mathcal
B} \setminus \{m\}$. We separate two cases, depending on the
number of double edges of the removed minor.

\begin{lemma} \label{quadrant}
Consider a weakened independence model obtained removing a
binomial with four double edges by a given family of adjacent
binomials ${\mathcal B}$, i.e. let ${\mathcal B}'={\mathcal
B}\setminus \lbrace m\rbrace$ where $m$ has four double edges. If
we let $A_{\mathcal{B}}$ be the orthogonal to $Z_{\mathcal B}$,
then the orthogonal to $Z_{\mathcal B'}$ is generated by the
elements of $A_{\mathcal{B}}$  and by the indicator vector $Q$ of
a quadrant centered on one of the indeterminates of the removed
minor.
\end{lemma}
\begin{proof}
First notice that the elements of $A_{\mathcal{B}}$ are orthogonal
to the columns of $Z_{\mathcal{B}'}$, i.e.
$A_{\mathcal{B}'}\supseteq A_{\mathcal{B}}$, and clearly, by Lemma
\ref{logveclemma}, one has
\begin{equation*}
\dim( A_{\mathcal{B}'})=\dim( A_{\mathcal{B}})+1 \, ,
\end{equation*}
where $A_{\mathcal{B}'}$ is the orthogonal to $Z_{\mathcal{B}'}$.
Now let $Q$ be the indicator vector of a quadrant centered on one
of the indeterminates of $m$. Then, $Q\not\in A_{\mathcal{B}}$ as
it is not orthogonal to the log-vector of $m$. But, $Q\in
A_{\mathcal{B}'}$ as each binomial in $\mathcal{B}'$ either avoid
the quadrant, or it is contained in the quadrant, or has exactly
two elements on the border of the quadrant. This is enough to
complete the proof.
\end{proof}

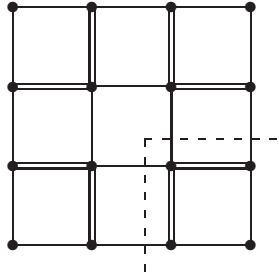
\begin{figure}
\begin{center}
\begin{picture}(110,110)(0,0)
\put(10,10){\circle*{4}} \put(10,10){\line(0,1){30}}
\put(10,10){\line(1,0){30}}

\put(10,40){\circle*{4}} \put(10,40){\line(0,1){30}}
\put(10,39){\line(1,0){30}} \put(10,41){\line(1,0){30}}

\put(10,70){\circle*{4}} \put(10,70){\line(0,1){30}}
\put(10,69){\line(1,0){30}} \put(10,71){\line(1,0){30}}

\put(10,100){\circle*{4}} \put(10,100){\line(1,0){30}}

\put(40,10){\circle*{4}} \put(39,10){\line(0,1){30}}
\put(41,10){\line(0,1){30}} \put(40,10){\line(1,0){30}}

\put(40,40){\circle*{4}} \put(40,40){\line(0,1){30}}
\put(40,40){\line(1,0){30}}

\put(40,70){\circle*{4}} \put(39,70){\line(0,1){30}}
\put(41,70){\line(0,1){30}} \put(40,70){\line(1,0){30}}

\put(40,100){\circle*{4}} \put(40,100){\line(1,0){30}}

\put(70,10){\circle*{4}}  \put(69,10){\line(0,1){30}}
\put(71,10){\line(0,1){30}} \put(70,10){\line(1,0){30}}

\put(70,40){\circle*{4}} \put(70,40){\line(0,1){30}}
\put(70,39){\line(1,0){30}} \put(70,41){\line(1,0){30}}

\put(70,70){\circle*{4}}  \put(69,70){\line(0,1){30}}
\put(71,70){\line(0,1){30}} \put(70,69){\line(1,0){30}}
\put(70,71){\line(1,0){30}}

\put(70,100){\circle*{4}} \put(70,100){\line(1,0){30}}

\put(100,10){\circle*{4}} \put(100,10){\line(0,1){30}}

\put(100,40){\circle*{4}} \put(100,40){\line(0,1){30}}

\put(100,70){\circle*{4}} \put(100,70){\line(0,1){30}}

\put(100,100){\circle*{4}}

\put(60,50){\dashline[8]{3}(0,0)(0,-50)}
\put(60,50){\dashline[8]{3}(0,0)(50,0)}
\end{picture}
\end{center}
\caption[]{Binomials of Example \ref{patexample}. The dashed line
delimits the quadrant defined in Lemma \ref{quadrant}.}
\label{quadfigure}
\end{figure}

The quadrant to be used in Example \ref{patexample} is sketched
Figure \ref{quadfigure}.

\begin{lemma} \label{external}
Consider a weakened independence model obtained removing a
binomial with not all the edges double  by a given family of
adjacent binomials ${\mathcal B}$, i.e. let ${\mathcal
B}'={\mathcal B}\setminus \lbrace m\rbrace$ where $m$ has not all
the edges double. If we let $A_{\mathcal{B}}$ be the orthogonal to
$Z_{\mathcal B}$, then the orthogonal to $Z_{\mathcal B'}$ is
generated by: the elements of $A_{\mathcal{B}}$, the indicator
vectors of the $MCC$s, of the $MCR$s and of the free cells.
\end{lemma}
\begin{proof}
Clearly, the elements of $A_{\mathcal{B}}$ are orthogonal to the
column of $Z_{\mathcal{B}'}$, i.e. $A_{\mathcal{B}'}\supseteq
A_{\mathcal{B}}$, and by Lemma \ref{logveclemma} one has
\begin{equation*}
\dim ( A_{\mathcal{B}'})=\dim (A_{\mathcal{B}})+1 \, ,
\end{equation*}
where $A_{\mathcal{B}'}$ is the orthogonal to $Z_{\mathcal{B}'}$.
The removed binomial $m$, with not all the edges double, can be
one of the following
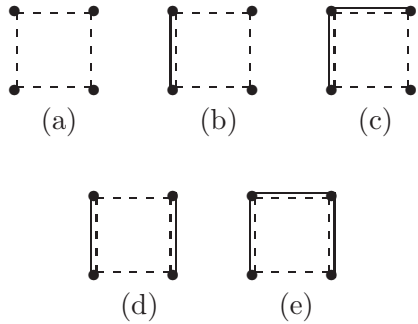
\begin{center}
\begin{picture}(170,135)(0,0)
\put(50,10){(d)}
\put(40,25){\circle*{4}}\put(40,55){\circle*{4}}\put(70,25){\circle*{4}}\put(70,55){\circle*{4}}
\put(41,26){\dashline[8]{3}(0,0)(0,28)(28,28)(28,0)(0,0)}
\put(39,24){\line(0,1){32}} \put(71,24){\line(0,1){32}}

\put(110,10){(e)}
\put(100,25){\circle*{4}}\put(100,55){\circle*{4}}\put(130,25){\circle*{4}}\put(130,55){\circle*{4}}
\put(101,26){\dashline[8]{3}(0,0)(0,28)(28,28)(28,0)(0,0)}
\put(99,24){\line(0,1){32}} \put(131,24){\line(0,1){32}}
\put(99,56){\line(1,0){32}}

\put(20,80){(a)}
\put(10,95){\circle*{4}}\put(10,125){\circle*{4}}\put(40,95){\circle*{4}}\put(40,125){\circle*{4}}
\put(11,96){\dashline[8]{3}(0,0)(0,28)(28,28)(28,0)(0,0)}

\put(80,80){(b)}
\put(70,95){\circle*{4}}\put(70,125){\circle*{4}}\put(100,95){\circle*{4}}\put(100,125){\circle*{4}}
\put(71,96){\dashline[8]{3}(0,0)(0,28)(28,28)(28,0)(0,0)}
\put(69,94){\line(0,1){32}}

\put(140,80){(c)}
\put(130,95){\circle*{4}}\put(130,125){\circle*{4}}\put(160,95){\circle*{4}}\put(160,125){\circle*{4}}
\put(131,96){\dashline[8]{3}(0,0)(0,28)(28,28)(28,0)(0,0)}
\put(129,94){\line(0,1){32}}\put(129,126){\line(1,0){32}}
\end{picture}
\end{center}
and to complete the proof we only need to present, in each case, a
vector $Q$ in $A_\mathcal{B'}$ which is not in $A_\mathcal{B}$. In
case (a), either we have a new free cell or not. If we have, let
$Q$ be the indicator vector of the free cell. Clearly, $Q\not\in
A_\mathcal{B}$, but $Q\in A_\mathcal{B'}$ has no minor is
involving the variable corresponding to the free cell. If we do
not have a new free cell, then we have a new $MMC$ or a new $MCR$
and its indicator vector is the required one. Repeating this kind
of argument in cases (b) through (e) we complete the proof.
\end{proof}

We are now ready to analyze the general case.

\begin{definition}
Let ${\mathcal B}$ be a set of adjacent minors and consider its
complement $\overline{\mathcal B}$ in the set of all adjacent
minors. Let $G_{\mathcal B}$ be the graph associated with
$\overline{\mathcal B}$. For each connected component of
$G_{\mathcal B}$ not touching the border of the table, we consider
the lex-smallest variable and we call it a {\em corner}.
\end{definition}

\begin{remark}
We notice that the number of corners is just the number of holes
one can find in the graph defined by the set of monomials.
\end{remark}

\begin{thm} \label{mainth}
Let ${\mathcal B}$ be a set of adjacent minors. Then the
orthogonal to $Z_{\mathcal B}$ is generated by the indicator
vectors of the $MCC$s, of the $MCR$s, of the free cells and by
quadrants centered in variables corresponding to corners.
\end{thm}
\begin{proof}
$\mathcal{B}$ can be constructed by the set of all the adjacent
minors by removing a minor at each time. Using Lemmas
\ref{quadrant} and \ref{external} we only need to show the result
for the set of all the adjacent minors, but this is a
straightforward consequence of Lemma \ref{logveclemma} and Theorem
\ref{maincounter}.
\end{proof}

\begin{remark}
In alternative to the straightforward use of Theorem \ref{mainth}
one can apply Lemmas \ref{external} and \ref{quadrant}  to
determine a sufficient statistic for a weakened independence model
with set of binomial ${\mathcal B}$. It is enough to start from
the complete set of adjacent minors and remove one by one the
minors not in ${\mathcal B}$. Notice that such an iterative
procedure, and the theorem itself, yield a system of generators of
the space orthogonal to $Z_{\mathcal B}$, but not a basis, i.e.
some of the vectors we add are redundant.
\end{remark}

We will show some examples and applications of Theorem
\ref{mainth} in the next sections.

\section{Exponential models} \label{expsect}

In Section \ref{suffsect} we have carried out some computations to
determine a sufficient statistic of a weakened independence model.
We are now able to find a parametric representation of the model.

Let us introduce unrestricted positive parameters $\zeta_1,
\ldots, \zeta_s$, where $s$ is the number of columns of the matrix
$A_{\mathcal B}$. If $A_{\mathcal B}$ has full rank, then $s$
coincide with the dimension of the vector sub-space orthogonal to
$Z_{\mathcal B}$.

The first step is to prove that weakened independence models
belong to the class of toric models and therefore they are
exponential (log-linear) models on the strictly positive simplex.
The first result in this direction is a rewriting of a general
theorem to be found in \cite{geiger|meek|sturmfels:06}.

\begin{remark}
The main result in \cite{geiger|meek|sturmfels:06} can be applied
to statistical models when the matrix representation $A_{\mathcal
B}$ of the sufficient statistic has non-negative entries.
Therefore, the theory developed in Section \ref{suffsect} is
relevant not only to actually determine a sufficient statistic,
but also to derive further theoretical properties of the weakened
independence models.
\end{remark}

Here we use again a vector notation in order to improve the
readability and simplify the formulae.

\begin{thm}[Geiger, Meek, Sturmfels (2006), Th. 3.2] \label{partheo}
Given a ${\mathcal B}$-weakened independence model, it can be
expressed as
\begin{equation} \label{parrepres}
p(\zeta) = \zeta^{A_{\mathcal B}}
\end{equation}
apart from the normalizing constant.
\end{thm}

Clearly, the parametrization in Eq. \eqref{parrepres} is not
unique. Theorem \ref{partheo} provides an easy way to switch from
the implicit representation to its parametrization. It is enough
to consider the matrix $A_{\mathcal B}$, whose columns are
orthogonal to the log-vectors of the binomials. As the columns of
$A_{\mathcal B}$ can be chosen with non-negative entries, then
each binomial in ${\mathcal B}$ vanishes in all points of the form
given in Eq. \eqref{parrepres}. This follows from a direct
substitution. The converse part is less intuitive and the proof is
not obvious.

As a corollary, Theorem \ref{partheo} allows us to consider
weakened independence models into the larger class of toric
models, as described in \cite{pistone|riccomagno|wynn:01} and
\cite{rapallo:07}. In the following result, we summarize the main
properties inherited from toric models.

\begin{proposition} \label{risriassunti}
Consider a ${\mathcal B}$-weakened independence model $V_{\mathcal
B}$.
\begin{enumerate}
\item $V_{\mathcal B}$ is a toric model;

\item With the constraint $p>0$, $V_{\mathcal B}$ is an
exponential model;

\item In case of sampling, the sufficient statistic for the sample
of size $n$ is the sum of the sufficient statistic of all
components of the sample.
\end{enumerate}
\end{proposition}
\begin{proof}
See Theorem 2 and the discussion in Section 3 of
\cite{rapallo:07}.
\end{proof}

\begin{remark}
As noticed in \cite{rapallo:07}, when we consider the general case
$p \geq 0$ instead of $p>0$, the toric model is not an exponential
model. Nevertheless it can be described as the disjoint union of a
suitable number of exponential models.
\end{remark}

We conclude this section with two examples.

\begin{example} \label{ex3per3}
As a first example, we consider a statistical model for $3 \times
3$ contingency tables defined through the binomials in Figure
\ref{fig3per3}.
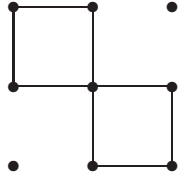
\begin{figure}
\begin{center}
\begin{picture}(80,80)(0,0)
\put(10,10){\circle*{4}}

\put(10,40){\circle*{4}} \put(10,40){\line(0,1){30}}
\put(10,40){\line(1,0){30}}

\put(10,70){\circle*{4}} \put(10,70){\line(1,0){30}}

\put(40,10){\circle*{4}} \put(40,10){\line(0,1){30}}
\put(40,10){\line(1,0){30}}

\put(40,40){\circle*{4}} \put(40,40){\line(0,1){30}}
\put(40,40){\line(1,0){30}}

\put(40,70){\circle*{4}}

\put(70,10){\circle*{4}}  \put(70,10){\line(0,1){30}}

\put(70,40){\circle*{4}}

\put(70,70){\circle*{4}}

\end{picture}
\end{center}
\caption[]{Binomials for Example \ref{ex3per3}.} \label{fig3per3}
\end{figure}
Using the theory developed in the previous section, the $MCR$s,
the $MCC$s and the free cells are sufficient to describe the
orthogonal $Z_{\mathcal B}$ and therefore the relevant matrices
are in Table \ref{matrici}.
\begin{table}
\begin{equation*}
[Z_{\mathcal B} \ | \ A_{\mathcal B}] =
\begin{pmatrix}
\begin{tabular}{cc|cccccccc}
1 & 0 & 1 & 0 & 0 & 1 & 0 & 0 & 0 & 0 \\
-1 & 0 & 1 & 0 & 0 & 0 & 1 & 0 & 0 & 0 \\
0 & 0 &  0 & 0 & 0 & 0 & 0 & 0 & 1 & 0 \\
-1 & 0 & 0  & 1 & 0 & 1 & 0 & 0 & 0 & 0 \\
1 & 1 & 0 & 1 & 0 & 0 & 1 & 0 & 0 & 0 \\
0 & -1 & 0 & 1 & 0 & 0 & 0 & 1 & 0 & 0 \\
0 & 0 & 0 & 0 & 0 & 0 & 0 & 0 & 0 & 1 \\
0  & -1 & 0 & 0 & 1 & 0 & 1 & 0 & 0 & 0 \\
0 & 0 & 0 & 0 & 1 & 0 & 0 & 1 & 0 & 0 \\
\end{tabular}
\end{pmatrix}
\end{equation*}
\caption{Matrices $Z_{\mathcal B}$ and $A_{\mathcal B}$ for
Example \ref{ex3per3}.}\label{matrici}
\end{table}
Thus, a parametrization with parameters $\zeta_1, \ldots, \zeta_8$
is:
\begin{equation*}
\left\{\begin{aligned}[l] p_{1,1} = & \zeta_1 \zeta_4 \\
p_{1,2} = & \zeta_1 \zeta_5 \\
p_{1,3} = & \zeta_7 \end{aligned} \right. \ \ \
\left\{\begin{aligned}[l] p_{2,1} = & \zeta_2 \zeta_4 \\
p_{2,2} = & \zeta_2 \zeta_5 \\
p_{2,3} = & \zeta_2 \zeta_6 \end{aligned} \right. \ \ \
\left\{\begin{aligned}[l] p_{3,1} = & \zeta_8 \\
p_{3,2} = & \zeta_3 \zeta_5 \\
p_{3,3} = & \zeta_3 \zeta_6 \end{aligned} \right.
\end{equation*}
\end{example}

\begin{example}
Let us consider the weakened independence model for $4 \times 4$
defined in Example \ref{patexample}. In that model, a minor with
four double edges has been removed and consequently the $MCR$s,
$MCC$s and free cells are not sufficient to describe the
orthogonal $Z_{\mathcal B}$. A vector must be added according to
Theorem \ref{quadrant}. Following the same approach as in the
previous example one can easily write down the matrices
$Z_{\mathcal B}$ and $A_{\mathcal B}$. A parametrization with
parameters $\zeta_1, \ldots, \zeta_9$ is:
\begin{equation*}
\left\{\begin{aligned}[l] p_{1,1} = & \zeta_1 \zeta_5 \\
p_{1,2} = & \zeta_1 \zeta_6 \\
p_{1,3} = & \zeta_1 \zeta_7 \\
p_{1,4} = & \zeta_1 \zeta_8 \end{aligned} \right. \ \ \
\left\{\begin{aligned}[l] p_{2,1} = & \zeta_2 \zeta_5 \\
p_{2,2} = & \zeta_2 \zeta_6 \\
p_{2,3} = & \zeta_2 \zeta_7 \\
p_{2,4} = & \zeta_2 \zeta_8 \end{aligned} \right. \ \ \
\left\{\begin{aligned}[l] p_{3,1} = & \zeta_3 \zeta_5 \\
p_{3,2} = & \zeta_3 \zeta_6 \\
p_{3,3} = & \zeta_3 \zeta_7 \zeta_9 \\
p_{3,4} = & \zeta_3 \zeta_8 \zeta_9 \end{aligned} \right. \ \ \
\left\{\begin{aligned}[l] p_{4,1} = & \zeta_4 \zeta_5 \\
p_{4,2} = & \zeta_4 \zeta_6 \\
p_{4,3} = & \zeta_4 \zeta_7 \zeta_9 \\
p_{4,4} = & \zeta_4 \zeta_8 \zeta_9 \end{aligned} \right.
\end{equation*}
\end{example}

\section{Inference and examples} \label{markovsect}

In the previous sections we have defined and studied the weakened
independence models. When the statistical model is given through a
set of adjacent minors ${\mathcal B}$, we are now able to compute
a sufficient statistic for a sample of size $1$ (see Proposition
\ref{risriassunti}) and to find a parametrization of the
statistical model $V_{\mathcal B}$ (see Theorem \ref{partheo}.) In
this section we give some ideas on how to compute maximum
likelihood estimates (MLE) and to perform exact inference through
Algebraic Statistics.

In the independence model defined by the set ${\mathcal C}$ of all
adjacent minors, the maximum likelihood estimate can be expressed
in closed form in terms of the observed value of the sufficient
statistic $T$. In the weakened independence models this is no
longer true, but numerical algorithms for log-linear models can be
used. As pointed out in Section \ref{expsect}, at least in the
strictly positive case, a weakened independence model is
log-linear and thus modified Newton-Raphson methods, Iterative
Proportional Fitting or EM methods can be used, see e.g.
\cite{agresti:02}. To compute the MLEs of of the cell
probabilities for the examples in this paper, we have used the
{\tt R} software, see \cite{rproject:06} , together with the
package {\tt gllm} (generalized log-linear models), see
\cite{duffy:06}. This package allows to define a generic matrix
for the sufficient statistic. This is the main advantage of the
{\tt gllm} package with respect to other available procedures in
different software systems.

Another theoretical result which highlights once again the
interplay between statistical models and polynomial algebra is the
Birch's Theorem (see e.g. \cite{bishop|fienberg|holland:75}). It
states that the MLE is the unique point $\hat p$ of the model
$V_{\mathcal B}$ which satisfies the constraints $A_{\mathcal B}^t
\hat p= A_{\mathcal B}^tp_{obs}$, where $p_{obs}$ are the observed
frequencies. We will have the opportunity to apply such result
later in Example \ref{swiss}.

Once the MLE is available, the goodness-of-fit can be evaluated
through a chi-squared test. The Pearson test statistic
\begin{equation}
C^2 = n \sum_{i,j} \frac {(p_{i,j}-\hat p_{i,j})^2} {\hat p_{i,j}}
\end{equation}
or the log-likelihood ratio test statistic
\begin{equation}
G^2 = 2 n \sum_{i,j} p_{i,j} \log \left( \frac {p_{i,j}} {{\hat
p}_{i,j}} \right)
\end{equation}
are evaluated and compared with the chi-square distribution with
$\#{\mathcal B}$ degrees of freedom, where $\#{\mathcal B}$ is the
cardinality of ${\mathcal B}$, see \cite{haberman:74} or
\cite{bishop|fienberg|holland:75}.

Alternatively, one can run the goodness-of-fit test within
Algebraic Statistics, using a Markov Chains Monte Carlo (MCMC)
algorithm. A number of papers have shown the relevance of this
approach, see \cite{diaconis|sturmfels:98}, and e.g.
\cite{rapallo:03}, \cite{aoki|takemura:05b}, and
\cite{chen|dinwoodie|dobra|huber:05}. The algebraic MCMC algorithm
was first described in \cite{diaconis|sturmfels:98}, and it is by
now widely used to compute non-asymptotic $p$-values for
goodness-of-fit tests in contingency tables problems.

Let $h$ be the observed contingency table for a sample of size
$n$, written as a vector in ${\mathbb N}^{IJ}$. The MCMC algorithm
is useful to efficiently sample from the reference set ${\mathcal
F}_t$ of a contingency table $h$ given a sufficient statistic with
matrix representation $A_{\mathcal B}$, i.e., from the set
\begin{equation*}
{\mathcal F}_t = \left\{h' \in {\mathbb N}^{IJ} \ | A_{\mathcal
B}^t h' = A_{\mathcal B}^t h \right\} \, .
\end{equation*}
The algorithm samples tables from ${\mathcal F}_t$ with the
appropriate hypergeometric distribution ${\mathcal H}$ through a
Markov chain based on a suitable set of moves making the chain
connected. Such set of moves is called a {\em Markov basis}. With
more details, a Markov basis is a set of tables $\{m_1, \ldots,
m_L\}$ with integer entries such that:
\begin{itemize}
\item $A_{\mathcal B}^t m_k = 0$ for all $k=1, \ldots, L$

\item if $h_1$ and $h_2$ are tables in ${\mathcal F}_t$, there
exist moves $m_{k_1}, \ldots, m_{k_A}$ and signs $\epsilon_1,
\ldots, \epsilon_A$ ($\epsilon_a= \pm1$) such that
\begin{equation*}
h_2 = h_1 + \sum_{a=1}^A \epsilon_a m_{k_a} \ \ \mbox{ and } \ \
h_1 + \sum_{a=1}^l \epsilon_a m_{k_a} \geq 0 \ \mbox{ for all }
l=1 , \ldots, A \, .
\end{equation*}
\end{itemize}
These conditions ensure the irreducibility of the Markov chain
defined by the algorithm below:
\begin{itemize}
\item Start from a table $h_1 \in {\mathcal F}_t$;

\item Choose a move $m_k$ uniformly in $\{m_1, \ldots, m_L\}$ and
a sign $\epsilon$ uniformly in $\{-1,1\}$. Define $h_2=h_1+
\epsilon m_k$;

\item Choose $u$ uniformly distributed in the interval $[0,1]$. If
$h_2 \geq 0$ and $\min\{1,H(h_2)/H(h_1)\}>u$ then move from $h_1$
to $h_2$, otherwise stay at $h_1$.
\end{itemize}

In the general case, the computation of a Markov basis needs
symbolic computations (the Diaconis-Sturmfels algorithm).
Nevertheless, a Markov basis for the weakened independence models
can be derived theoretically. In the following, we will determine
a Markov basis for the weakened independence models.

Given the set ${\mathcal B}$ of binomials defining the ${\mathcal
B}$-weakened independence model, let $A_{\mathcal B}$ be the
matrix representation of the sufficient statistic. Moreover, we
denote by ${\mathcal I}_{\mathcal B}$ the polynomial ideal in
${\mathbb R}[p]$ generated by the binomials in ${\mathcal B}$.
Diaconis and Sturmfels (\cite{diaconis|sturmfels:98}, Theorem
$3.1$) proved that a Markov basis is formed by the log-vectors of
a set of generators of the toric ideal associated to $A_{\mathcal
B}$, i.e. the ideal
\begin{equation*}
{\mathcal J}_{\mathcal B} = \{p^a-p^b \ | \ a,b \in {\mathbb
R}^{IJ} \ , \ A_{\mathcal B}^t(a) = A_{\mathcal B}^t(b) \}.
\end{equation*}

Therefore, the computation of a Markov basis translates into the
computation of a set of generators of a toric ideal.

Bigatti {\it et al.} \cite{bigatti|lascala|robbiano:99} showed
that the toric ideal associated to $A_{\mathcal B}$ is the
saturation of ${\mathcal I}_{\mathcal B}$ with respect to the
product of the indeterminates. Such ideal is defined as:
\begin{equation*}
{\mathcal I}_{\mathcal B}:(p_{1,1}\cdots
p_{I,J})^{\infty}=\left\lbrace f \in\mathbb{R}[p] \ | \
(p_{1,1}\cdots p_{I,J})^n f \in {\mathcal I}_{\mathcal B} \mbox{
for some } n \right\rbrace.
\end{equation*}

In order to compute a set of generators of this ideal, one can use
symbolic algebra packages, e.g. the function {\tt Toric} of CoCoA.
For further details on ideals and their operations, see e.g.
\cite{cox|little|oshea:92} and \cite{kreuzer|robbiano:00}.

\begin{example}
The data we present as a first example in this section have been
collected by one of the authors in his Biostatistics course. Each
of the $34$ students must submit a homework before the exam and
this report is evaluated by two instructors on a scale with levels
$\{1,2,3\}$. The final grade is the maximum of the two
evaluations. The data are in the Table \ref{biostat}.
\begin{table}
\begin{center}
\begin{tabular}{ccccc|c}
\hline\hline
{Second instructor}& \ \  & \multicolumn{3}{c}{First instructor} & \\
\hline  & & $1$ & $2$ & $3$ & Total \\ \hline
  $1$ & & $7$ & $5$ & $0$ & $12$ \\
    & & $(6.52)$ & $(5.48)$ & $(0)$ &  \\
 $2$ & & $4$ & $5$ & $2$ & $11$ \\
     & & $(4.48)$ & $(3.76)$ & $(2.76)$ &  \\
  $3$ & & $1$ & $5$ & $5$ & $11$ \\
      & & $(1)$ & $(5.76)$ & $(4.24)$ & \\ \hline
 Total & & 12 & 15 & 7 & 34
\end{tabular}
\end{center}
\caption{Evaluation of $34$ homeworks by two instructors. In
parentheses are the the MLE estimates for the weakened
independence model.}\label{biostat}
\end{table}

The model we use to analyze such data is the model defined by the
adjacent minors in Example \ref{ex3per3}. Using the {\tt gllm}
package, we obtain the MLEs written in parentheses in the table.
The Pearson statistic is $0.9863$. Running a MCMC algorithm with a
Markov basis consisting of $2$ moves, we find a $p$-value of
$0.6665$ for the goodness-of-fit test, showing a good fit. The
Monte Carlo computations are based on a sample of $10,000$ tables,
with a burn-in phase of $50,000$ tables and sampling every $50$
steps.

Models of this kind are used in \cite{carlini|rapallo:08} to
detect category indistinguishability both in intra-rater and in
inter-rater agreement problems. The model we used shows that
categories $1$ and $2$ are confused, as well as categories $2$ and
$3$. This lack of distinguishability can be ascribed to a relevant
non-homogeneity of the marginal distributions.
\end{example}

\begin{example}
The data in Table \ref{cholest} show the cross-classification of
$103$ subjects with respect to $2$ ordinal variables: the smoking
level, $4$ categories from ``No Smoking'' to ``More than $10$
cigarettes'', and the quantity of High-Density Lipoprotein
Cholesterol (HDLP) in the blood, $4$ categories from ``Normal'' to
``Abnormal''. The data are presented in \cite{jeong|jhun|kim:05}
and analyzed by the authors under both the independence model and
the RC (Row Column effects) model. The authors compute the exact
$p$-values for the independence model (0.049) and for the RC model
(0.657) using the log-likelihood ratio test statistic.
\begin{table}
\begin{center}
\begin{tabular}{cccccc|c} \hline\hline
Smoking level & \ \ & \multicolumn{4}{c}{HDLC} & \\ \hline
 & & $1$ & $2$ & $3$ & $4$ & Total \\ \hline
 $1$  & & $15$ & $3$ & $6$ & $1$ & 25 \\
 $2$  & & $8$ & $4$ & $7$ & $2$  & $21$ \\
 $3$  & & $11$ & $6$ & $15$ & $3$  & $35$ \\
 $4$ & & $5$ & $1$ & $11$ & $5$  & $22$ \\ \hline
  Total   &  &  $39$ & $14$ & $39$ & $11$ & $103$
\end{tabular}
\end{center}
\caption{Cross-classification of Smoking level and HDLC. Smoking
levels: $1=$``No Smoking'', $2=$``Less than $5$ cigarettes'',
$3=$``Less than $10$ cigarettes'', $4=$``More than $10$
cigarettes''. HDLC levels: $1=$``Normal'', $2=$``Low Normal'',
$3=$``Borderline'', $4=$``Abnormal''.}\label{cholest}
\end{table}

We use a weakened independence model with binomials in Figure
\ref{figcholest}. According to Theorem \ref{mainth}, a sufficient
statistic is formed by $4$ $MCR$s, $4$ $MCC$s and $3$ free cells.
Moreover, the relevant Markov basis has $14$ binomials.

\begin{figure}
\begin{center}
\begin{picture}(110,110)(0,0)
\put(10,10){\circle*{4}}\put(10,10){\line(0,1){30}}
\put(10,10){\line(1,0){30}}

\put(10,40){\circle*{4}}\put(10,40){\line(0,1){30}}
\put(10,39){\line(1,0){30}} \put(10,41){\line(1,0){30}}

\put(10,70){\circle*{4}} \put(10,70){\line(1,0){30}}

\put(10,100){\circle*{4}}

\put(40,10){\circle*{4}} \put(39,10){\line(0,1){30}}
\put(41,10){\line(0,1){30}} \put(40,10){\line(1,0){30}}

\put(40,40){\circle*{4}} \put(39,40){\line(0,1){30}}
\put(41,40){\line(0,1){30}} \put(40,39){\line(1,0){30}}
\put(40,41){\line(1,0){30}}

\put(40,70){\circle*{4}} \put(40,70){\line(0,1){30}}
\put(40,69){\line(1,0){30}} \put(40,71){\line(1,0){30}}

\put(40,100){\circle*{4}} \put(40,100){\line(1,0){30}}

\put(70,10){\circle*{4}} \put(70,10){\line(0,1){30}}

\put(70,40){\circle*{4}} \put(70,40){\line(0,1){30}}

\put(70,70){\circle*{4}} \put(69,70){\line(0,1){30}}
\put(71,70){\line(0,1){30}} \put(70,70){\line(1,0){30}}

\put(70,100){\circle*{4}} \put(70,100){\line(1,0){30}}

\put(100,10){\circle*{4}}

\put(100,40){\circle*{4}}

\put(100,70){\circle*{4}} \put(100,70){\line(0,1){30}}

\put(100,100){\circle*{4}}
\end{picture}
\end{center}
\caption[]{Weakened independence model for the cholesterol data.}
\label{figcholest}
\end{figure}
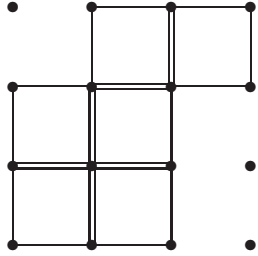
With our model we find a $p$-value of $0.7205$. Therefore, this
weakened independence model fits better than the independence
model and it is as good as the RC model. In particular, removing
only three adjacent minors from the complete configuration with 9
minors, we obtain a model whose fit is dramatically improved.
Moreover, the removed minors allow us to identify quickly the
cells which cause the departure from independence.
\end{example}

\begin{example}\label{swiss}
To conclude this section, we show how the models defined in the
present paper have some relationships with a recent problem, first
stated by Bernd Sturmfels in $2005$ and known as the ``100 Swiss
Francs Problem'', see \cite{sturmfels:07}. Such problem is related
to the modelling of DNA sequence alignments. We briefly describe
the probabilistic experiment. For a plain description, the reader
can refer to \cite{patcher|sturmfels:05}, Example $1.15$.

A DNA sequence is a sequence of symbols in the alphabet $\{${\tt
A}, {\tt T}, {\tt C}, {\tt G}$\}$. A major problem in molecular
biology is to compare two DNA sequences. In \cite{sturmfels:07},
the following observed sequences were considered:
\begin{center}
{\tt ATCACCAAACATTGGGATGCCTGTGCATTTGCAAGCGGCT} \\
{\tt ATGAGTCTTAAAACGCTGGCCATGTCCATCTTAGACAGCG}
\end{center}
leading to the observed table below:
\begin{equation} \label{DNAdata}
\begin{pmatrix}
\begin{tabular}{cccc}
4 & 2 & 2 & 2 \\
2 & 4 & 2 & 2 \\
2 & 2 & 4 & 2 \\
2 & 2 & 2 & 4
\end{tabular}
\end{pmatrix}
\end{equation}
The hypothesis of the author is that such two DNA sequences are
generated through a (biased) coin and four tetrahedral dice $D_1,
D_2, D_3, D_4$ with the letters {\tt A}, {\tt T}, {\tt C}, {\tt G}
on the facets. When the coin outcome is ``Head'', then the dice
$D_1$ and $D_2$ are rolled. The outcome of $D_1$ is registered in
the first sequence and the outcome of $D_2$ in the second one.
When the coin outcome is ``Tail'', then the dice $D_3$ and $D_4$
are rolled.

Denote by $q_1, \ldots , q_4$ the probability vectors of the four
dice and by $\alpha$ the probability of ``Head'' in the coin.
Then, the probability distribution of the final outcome is
\begin{equation}
\alpha q_1q_2^t + (1 - \alpha)q_3q_4^t \, .
\end{equation}

Therefore, the construction of this experiment leads us to
consider the statistical model  of $4 \times 4$ matrices of
probabilities whose rank is less than or equal to $2$. The author
conjectured that the maximum likelihood estimate of the
probabilities under the model of matrices with rank at most $2$
is:
\begin{equation} \label{DNAmaxlike}
\hat P_g = \frac {1} {40}\begin{pmatrix}
\begin{tabular}{cccc}
3 & 3 & 2 & 2 \\
3 & 3 & 2 & 2 \\
2 & 2 & 3 & 3 \\
2 & 2 & 3 & 3
\end{tabular}
\end{pmatrix}
\end{equation}
Further analyses to be found in
\cite{fienberg|hersh|rinaldo|zhou:07} show that the model is
non-identifiable and that numerical methods are able to identify
$3$ global maxima and $4$ local maxima for the likelihood
function. Apart from the simultaneous permutation of the row and
column labels, the global maximum is reached for the matrix in Eq.
\eqref{DNAmaxlike} and the local maximum is obtained by the
matrix:
\begin{equation} \label{DNAmaxlocal}
\hat P_l =\frac {1} {40}\begin{pmatrix}
\begin{tabular}{cccc}
8/3 & 8/3 & 8/3 & 2 \\
8/3 & 8/3 & 8/3 & 2 \\
8/3 & 8/3 & 8/3 & 2 \\
2 & 2 & 2 & 4
\end{tabular}
\end{pmatrix}
\end{equation}

Now, consider the following probability models for $4 \times 4$
contingency tables:
\begin{itemize}
\item the model $M$ of the matrices with rank at most $2$;

\item the weakened independence models $M_1$ and $M_2$ whose
binomials are presented in Figure \ref{fig100francs}.
\end{itemize}

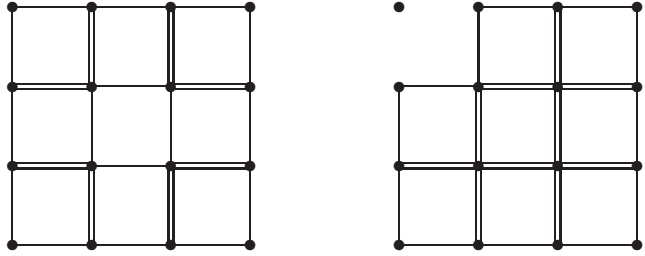
\begin{figure}
\begin{center}
\begin{tabular}{ccc}
\begin{picture}(110,110)(0,0)
\put(10,10){\circle*{4}} \put(10,10){\line(0,1){30}}
\put(10,10){\line(1,0){30}}

\put(10,40){\circle*{4}} \put(10,40){\line(0,1){30}}
\put(10,39){\line(1,0){30}} \put(10,41){\line(1,0){30}}

\put(10,70){\circle*{4}} \put(10,70){\line(0,1){30}}
\put(10,69){\line(1,0){30}} \put(10,71){\line(1,0){30}}

\put(10,100){\circle*{4}} \put(10,100){\line(1,0){30}}

\put(40,10){\circle*{4}} \put(39,10){\line(0,1){30}}
\put(41,10){\line(0,1){30}} \put(40,10){\line(1,0){30}}

\put(40,40){\circle*{4}} \put(40,40){\line(0,1){30}}
\put(40,40){\line(1,0){30}}

\put(40,70){\circle*{4}} \put(39,70){\line(0,1){30}}
\put(41,70){\line(0,1){30}} \put(40,70){\line(1,0){30}}

\put(40,100){\circle*{4}} \put(40,100){\line(1,0){30}}

\put(70,10){\circle*{4}}  \put(69,10){\line(0,1){30}}
\put(71,10){\line(0,1){30}} \put(70,10){\line(1,0){30}}

\put(70,40){\circle*{4}} \put(70,40){\line(0,1){30}}
\put(70,39){\line(1,0){30}} \put(70,41){\line(1,0){30}}

\put(70,70){\circle*{4}}  \put(69,70){\line(0,1){30}}
\put(71,70){\line(0,1){30}} \put(70,69){\line(1,0){30}}
\put(70,71){\line(1,0){30}}

\put(70,100){\circle*{4}} \put(70,100){\line(1,0){30}}

\put(100,10){\circle*{4}} \put(100,10){\line(0,1){30}}

\put(100,40){\circle*{4}} \put(100,40){\line(0,1){30}}

\put(100,70){\circle*{4}} \put(100,70){\line(0,1){30}}

\put(100,100){\circle*{4}}
\end{picture}

& \ \ \ \ \ &

\begin{picture}(110,110)(0,0)
\put(10,10){\circle*{4}}\put(10,10){\line(0,1){30}}
\put(10,10){\line(1,0){30}}

\put(10,40){\circle*{4}}\put(10,40){\line(0,1){30}}
\put(10,39){\line(1,0){30}} \put(10,41){\line(1,0){30}}

\put(10,70){\circle*{4}} \put(10,70){\line(1,0){30}}

\put(10,100){\circle*{4}}

\put(40,10){\circle*{4}} \put(39,10){\line(0,1){30}}
\put(41,10){\line(0,1){30}} \put(40,10){\line(1,0){30}}

\put(40,40){\circle*{4}} \put(39,40){\line(0,1){30}}
\put(41,40){\line(0,1){30}} \put(40,39){\line(1,0){30}}
\put(40,41){\line(1,0){30}}

\put(40,70){\circle*{4}} \put(40,70){\line(0,1){30}}
\put(40,69){\line(1,0){30}} \put(40,71){\line(1,0){30}}

\put(40,100){\circle*{4}} \put(40,100){\line(1,0){30}}

\put(70,10){\circle*{4}} \put(69,10){\line(0,1){30}}
\put(71,10){\line(0,1){30}} \put(70,10){\line(1,0){30}}

\put(70,40){\circle*{4}} \put(69,40){\line(0,1){30}}
\put(71,40){\line(0,1){30}} \put(70,39){\line(1,0){30}}
\put(70,41){\line(1,0){30}}

\put(70,70){\circle*{4}} \put(69,70){\line(0,1){30}}
\put(71,70){\line(0,1){30}} \put(70,69){\line(1,0){30}}
\put(70,71){\line(1,0){30}}

\put(70,100){\circle*{4}} \put(70,100){\line(1,0){30}}

\put(100,10){\circle*{4}} \put(100,10){\line(0,1){30}}

\put(100,40){\circle*{4}} \put(100,40){\line(0,1){30}}

\put(100,70){\circle*{4}} \put(100,70){\line(0,1){30}}

\put(100,100){\circle*{4}}
\end{picture}
\end{tabular}
\end{center}
\caption[]{Weakened independence models $M_1$ (left) and $M_2$
(right) for the $100$ Swiss francs problem.}\label{fig100francs}
\end{figure}

One can easily check that both $M_1$ and $M_2$ are proper subsets
of $M$. In fact, in $M_1$ the first row is proportional to the
second row and so are the third and the fourth, while in $M_2$,
the first three rows are proportional.

Now it is easy to check by direct substitution that the matrix
$\hat P_g$ of global maximum in Eq. \eqref{DNAmaxlike} belongs to
$M_1$, while the matrix $\hat P_l$ of local maximum in Eq.
\eqref{DNAmaxlocal} belongs to $M_2$. Such matrices are the MLEs
for the two models, respectively.
\end{example}

\begin{remark}
This is not a proof that the matrix $\hat P_g$ is the MLE for the
model $M$. Nevertheless, it is interesting to notice that our
models contain both the local and the global maxima. However, our
model can not suffice to find the local extrema in $M$ as $M$ is
not a toric model, e.g., it is not defined by binomials.
\end{remark}

\section{Final remarks and future work} \label{finremsect}

In this paper we used the binomial representation of the
independence model to introduce a new class of statistical
models:these models are devised to weaken independence. We studied
their sufficient statistic, we proved that they are log-linear
models for strictly positive probabilities and we showed how to
make inference on these models. Some numerical examples emphasized
the importance and the wide applicability of our models.

We have in mind different ways to generalize this work. First, we
want to find a procedure to characterize the relevant Markov bases
to be used in the Diaconis-Sturmfels algorithm. Then, we plan to
consider models defined by non-adjacent $2 \times 2$ minors and
try to analyze them with similar techniques. Moreover, we are
interested in the study of higher dimensional minors, e.g., $3
\times 3$ minors which appear in the definition of the models in
Example \ref{swiss}. Finally, for large tables there will be many
weakened independence models and model selection strategies must
be studied. Further applications of this kind of models will be
investigated, in particular in the field of computational biology.
From a geometrical point of view, we would like to explore the
structure of the varieties defined by some adjacent minors, as
done in \cite{hosten|sullivant:04} when all adjacent minors are
considered.

\bibliographystyle{acmtrans-ims}
\bibliography{tutto2}

\end{document}